\documentclass[12pt]{amsart}
\usepackage{amssymb}
\ifx\mathrm\undefined\let\mathrm\rm\fi
\ifx\mathbf\undefined\let\mathbf\bf\fi
\ifx\mathfrak\undefined\let\mathfrak\frak\fi
\ifx\mathcal\undefined\let\mathcal\cal\fi
\ifx\mathbb\undefined\let\mathbb\Bbb\fi
\ifx\emph\undefined\let\emph\it\fi

\newcommand{\Z}{{\mathbb Z}}

\newcommand{\R}{{\mathbb R}}
\newcommand{\C}{{\mathbb C}}
\newcommand{\Q}{{\mathbb Q}}

\newcommand{\Ref}[1]{{(\ref{#1})}}
\newcommand{\be}{\begin{displaymath}}
\newcommand{\ee}{\end{displaymath}}
\newcommand{\bea}{\begin{eqnarray*}}
\newcommand{\eea}{\end{eqnarray*}}

\newcommand{\bean}{\begin{eqnarray}}
\newcommand{\eean}{\end{eqnarray}}

\font\bb=msbm10 at9.98pt
\def\semidirect{\hbox{$\;$\bb\char'156$\;$}}

\newcommand{\dontprint}[1]{\relax}

\newenvironment{prf}{\noindent{\it Proof\/}:}{$\;\square$
\par\medskip}

\newtheorem%
{thm}{Theorem}[section]
\newtheorem%
{proposition}[thm]{Proposition}
\newtheorem%
{lemma}[thm]{Lemma}

\newtheorem%
{lemmadef}[thm]{Lemma-Definition}
\newtheorem%
{corollary}[thm]{Corollary}
\newtheorem%
{conjecture}[thm]{Conjecture}

\title[The Multiplication Formula for the Elliptic Gamma Function]
{Multiplication Formulas for the Elliptic Gamma Function}
\author[G. Felder and A. Varchenko]
{Giovanni Felder${}^{*}$ 
\and Alexander Varchenko${}^{**,1}$}
\thanks{${}^1$Supported in part by NSF grant  DMS-9801582}
\date{November 2002}
\begin{document}
\maketitle
\centerline{\it ${}^*$Departement Mathematik, ETH-Zentrum,}
\centerline{\it 8092
Z\"urich, Switzerland}
\centerline{felder@math.ethz.ch}
\medskip
\centerline{\it ${}^{**}$Department of Mathematics,
University of North Carolina at Chapel Hill,}
\centerline{\it Chapel Hill, NC 27599-3250, USA}
\centerline{anv@email.unc.edu}
\newcommand{\sig}{\sigma}
\begin{abstract}
The elliptic gamma function is a generalization of the Euler gamma
function.  Its trigonometric
and rational degenerations are the Jackson q-gamma function and the
Euler gamma function.  We prove 
multiplication formulas for the elliptic gamma function, whose
degenerations are the Gauss-Askey multiplication formula for the Euler and trigonometric
gamma functions.
\end{abstract}
%\section{Introduction}

\bigskip

\section{Introduction} Special functions defined by infinite products often
have duplication formulas.
Here are some examples.
\bean
\text{sin}\, (2\pi z)\ &=&\ 2\ 
\text{sin}\, (\pi z)\ \text{sin} \, (\pi ( z + \frac{1}{2}))\ ,
\notag
\\
\Gamma ( 2z) \ \sqrt{\pi} \ & =&\
2^{2z-1}\ \Gamma(z)\ \Gamma ( z + \frac{1}{2})\ ,
\\
\Gamma_q ( 2z) \ \Gamma_{q^2}(\frac{1}{2}) \ & = & \
[2]_q^{2z-1}\ \Gamma_{q^2}(z)\ \Gamma_{q^2} ( z + \frac{1}{2})\ ,
\\
\theta_0(2z, \tau)\ &= & \ \theta_0( z, \tau)\ \theta_0(z + \frac{1}{2},\tau ) \
\theta_0(z + \frac{\tau}{2},\tau) \ \theta_0(z + \frac{1 + \tau}{2},\tau) \ .
\eean

The function $\Gamma(z) = \int_0^\infty
t^{z-1}e^{-t}dt$  is the Euler gamma function. It
 satisfies the functional equation $\Gamma (z+1)\ =\
z\ \Gamma(z)$. Formula (1) is Legendre's duplication formula.

The function $\Gamma_q (z)$ is Jackson's $q$-gamma function. Set $x = e^{2\pi i\,z}$, 
 $q=e^{2\pi i \tau}$, and 
denote 
\[
(x;q)\ =\ \prod_{j=0}^\infty(1-xq^j) .
\]
Then
\[
\Gamma_q (z) \ =\ \Gamma_{\mathrm{trig}}(z,\tau)\ = \ 
(1-q)^{1-z}\ \frac{(q;q)}{(q^z;q)}\ .
\]
The $q$-gamma function obeys the functional equation
\[
\Gamma_q(z+1)\ =\ [z]_q \ 
\Gamma_{q}(z) \ ,
\]
where $[z]_q = \frac {1-e^{2\pi i\tau\, z}} {1-e^{2\pi i\tau}} $
is the trigonometric analog of the number $z$. The $q$-gamma function
degenerates to Euler's gamma function,
\[
\lim_{\tau\to 0}\ \Gamma_{\mathrm{trig}}(z,\tau)\
=\ \Gamma(z)\ .
\]
Formula (2) is Askey's duplication formula, see \cite{A}.

The function $\theta_0(z, \tau) \ =\ (x, q) (q/x , q)$ in (3) 
is one of Jacobi's theta functions. Formula (3) see for instance in \cite{Ra}.

In this paper we give two duplication formulas for elliptic analogs of
the gamma function, 
\bea 
 \Gamma(2z,\tau, \sig) = && 
\Gamma(z, \tau,\sig)
 \Gamma(z+\frac{\tau}{2}, \tau, \sig) \Gamma(z+\frac{\sig}{2},
\tau, \sig) \Gamma(z+\frac{\tau + \sig}{2}, \tau, \sig) 
\\ 
&&
\Gamma(z+\frac{1}{2}, \tau, \sig) \Gamma(z+\frac{1+\tau}{2}, \tau, \sig)
\Gamma(z+\frac{1+\sig}{2}, \tau, \sig) \Gamma(z+\frac{1+\tau + \sig}{2},
\tau, \sig)\ , 
\eea
\bea
\bar\Gamma(2z,\tau, \sig) \bar\Gamma(\frac{1}{2}, 2\tau,  \sig)
 =   
\Bigg( \frac{\theta_0(2 \tau, \sig)}{\theta_0(\tau,\sig)} 
\Bigg)^{2z-1} \bar\Gamma(z, 2\tau, \sig) 
\bar\Gamma(z+\frac{1}{2}, 2\tau, \sig)\ , 
\eea 
see definitions below. The expression $\frac{\theta_0( z \tau
,\sig)}{\theta_0(\tau,\sig)}$ is an elliptic analog of the number $z$.
We have the trigonometric and rational limits of the
theta function:
\[
\frac{\theta_0(z \tau ,\sig)}{\theta_0(\tau,\sig)}
\ 
\stackrel{\sig \to i\infty}{\longrightarrow}
\
\frac
{1-e^{2\pi i\tau\, z}}
{1-e^{2\pi i\tau}}
\
\stackrel{\tau \to 0}{\longrightarrow}
\ z\ .
\]

\section{Elliptic gamma function}
The {\em elliptic
gamma function} is an elliptic generalization
of the Euler gamma function. It is the meromorphic function
of three complex variables $z,\tau,\sig$, with 
$\mathrm{Im}\,\tau,\mathrm{Im}\,\sig>0$ defined
by the convergent infinite product
\[
\Gamma(z,\tau,\sig)=\prod_{j,k=0}^\infty
\frac
{1-e^{2\pi i((j+1)\tau+(k+1)\sig-z)}}
{1-e^{2\pi i(j\tau+k\sig+z)}}\,.
\]
It is the unique solution
of a functional equation involving the Jacobi
theta function $\theta_0$. 
\bigskip

{\bf Theorem.} \ \cite{FV1}\
{\it Suppose that $\tau,\sig$ are complex numbers with positive 
imaginary part. Then $u(z)=\Gamma(z,\tau,\sig)$ is the
unique meromorphic solution of the difference equation
\[
u(z+\sig)=\theta_0(z,\tau)u(z)
\]
such that:
\begin{enumerate}
\item[(i)] $u(z)$ obeys $u(z+1)=u(z)$ and is holomorphic
on the upper half plane
$\mathrm{Im}\, z>0$,
\item[(ii)] $u((\tau+\sig)/2)=1$. 
\end{enumerate}
}

\bigskip

The elliptic gamma function first appeared in  \cite{R}.
The modular properties of the elliptic gamma function and their
relations to $SL(3, \Z)$ are discussed in \cite {FV1}, appearances and applications 
of the elliptic gamma function can be found in \cite{B, DP, Kyoto1, Kyoto2, FTV, FV2,
FV3, FV4}.

\bigskip

Let
$\bar\Gamma$ be the function
\[
\bar\Gamma(z,\tau,\sig)\ =\ \frac{(q;q)}{(r;r)}\
\theta_0(\tau, \sig)^{1-z}\
\Gamma(z \tau,\tau,\sig)\ ,
\qquad 
q = e^{2\pi i\tau},\
r = e^{2\pi i\sig}.
\]
Then $u(z)=\bar\Gamma(z,\tau,\sig)$ is a solution of the functional equation
\[
u(z+1)\ =\ \frac{\theta_0(\tau z,\sig)}{\theta_0(\tau,\sig)}\ u(z).
\]
The normalization was chosen here so that $u(1)=1$.
As $\sig \to i\infty$ we recover  Jackson's  $q$-gamma function,
\[
\Gamma_{\mathrm{trig}}(z,\tau)\ =\
\lim_{\sig\to i\infty}\ 
\bar\Gamma(z,\tau,\sig)\ .
\]

\bigskip

\section{Multiplication Formulas}

\subsection{The first multiplication formula}
${ }$ \newline ${ }$

{\bf Theorem.} {\it For any natural $n$ we have}
\[
\Gamma(nz,\tau, \sig)\ =
\
\prod_{k_1,k_2,k_3\,=\,0}^{n-1} \
\Gamma(z+\frac{k_1+k_2\tau+k_3\sig}{n}, \tau, \sig)\ .
\]

{\bf Proof.} Let $w = e^{2\pi i/n}$. Then the right hand side of this formula is
\bea
&&
\prod_{k_1,k_2,k_3\,=\,0}^{n-1} \ 
\prod_{l,m = 0}^\infty 
\frac 
{1-  w^{-k_1} q^{l +  \frac{n-k_2}{n}} r^{m + \frac{n-k_3}{n}}
x^{-1}}
{1-  w^{k_1} q^{l +  \frac{k_2}{n}} r^{m + \frac{k_3}{n}}x } 
\\
&&
=
\prod_{k_2,k_3\,=\,0}^{n-1} \ 
\prod_{l,m = 0}^\infty 
\frac 
{1-  q^{(l+1)n - k_2 } r^{(m+1)n - k_3 } x^{-n}}
{1-  q^{ln + k_2} r^{mn + k_3}x^{n}} =
\prod_{l,m = 0}^\infty 
\frac 
{1-  q^{l+1} r^{m+1} x^{-n}}
{1-  q^{l} r^{m}x^{n}} 
=
\Gamma (nz, \tau, \sig) .
\eea
$\square$

\bigskip
\subsection{The second  multiplication formula} 
${ }$ \newline ${ }$

{\bf Theorem.}\ {\it For any natural $n$ we have}
\bea
&&
\bar\Gamma(nz,\tau, \sig)\ \bar\Gamma(\frac{1}{n}, n\tau, \sig)
\ \bar\Gamma(\frac{2}{n}, n\tau,  \sig)\ \dots 
\ \bar\Gamma(\frac{n-1}{n}, n\tau, \sig)\ =
\\
&&
\Bigg(
\frac{\theta_0(n\tau ,\sig)}{\theta_0(\tau,\sig)}
\Bigg)^{nz-1}\
\bar\Gamma(z, n\tau, \sig)\
\bar\Gamma(z+\frac{1}{n}, n\tau, \sig)
\ \bar\Gamma(z+\frac{2}{n}, n\tau, \sig)\ \dots 
\ \bar\Gamma(z+\frac{n-1}{n}, n\tau, \sig)\ .
\eea

\bigskip

The theorem is an easy consequence of the following two lemmas.

\bigskip

{\bf Lemma.}\
{\it For any natural $m$ and $ n$, we have}
\[
\Gamma(z, \tau, \sig)\ = \
\prod_{a=0}^{m-1}
\prod_{b=0}^{n-1}\
\Gamma(z + a\tau + b\sig, m\tau, n\sig)\  .
\]

\bigskip

{\bf Lemma.} \ {\it
For any natural $ n$,  we have}
\[
\Gamma(\tau, n\tau, \sig)\ 
\Gamma(2\tau, n\tau, \sig)\ \dots
\ \Gamma((n-1)\tau, n\tau, \sig)\ 
=\
\frac{1}{(q, q^n) (q^2, q^n) \dots (q^{n-1}, q^n)}\ .
\]

\end{document}